\documentclass[twoside,12pt]{article}
\usepackage{amsmath,amssymb,epsfig}
\date{}
\newtheorem{Lemma}{LEMMA}[section]
\newtheorem{Corollary}[Lemma]{COROLLARY}
\newtheorem{Theorem}[Lemma]{THEOREM}
\newtheorem{Proposition}[Lemma]{PROPOSITION}
\newtheorem{Definition}[Lemma]{DEFINITION}

\newcommand{\bnum}{\begin{enumerate}}
\newcommand{\enum}{\end{enumerate}}
\newcommand{\bi}{\begin{itemize}}
\newcommand{\ei}{\end{itemize}}
\newcommand{\btab}{\begin{tabular}}
\newcommand{\etab}{\end{tabular}}
\newcommand{\beq}{\begin{eqnarray*}}
\newcommand{\eeq}{\end{eqnarray*}}
\newcommand{\beqn}{\begin{eqnarray}}
\newcommand{\eeqn}{\end{eqnarray}}

\newcommand{\bq}{\begin{equation}}
\newcommand{\eq}{\end{equation}}

\newcommand{\CC}{{\cal C}}

\newcommand{\CL}{{\cal L}}

\newcommand{\CT}{{\cal T}}

\def\phi{\varphi}
\def\epsilon{\varepsilon}

\newcommand{\BR}{\mathbb R}

\newcommand{\kasten}{\vbox{\hrule height 8pt width 8.6pt depth -7.4pt
    \hbox{\vrule width 0.6pt height 7.4pt
    \kern 7.4pt \vrule width 0.6pt height 7.4pt}
    \hrule height 0.6pt width 8.6pt}}
\newcommand{\ok}{\hfill\kasten}
\newcommand{\bpf}{\begin{Proof}}
\newcommand{\epf}{\ok\end{Proof}\bigskip\noindent}
\newcommand{\bthm}{\begin{Theorem}}
\newcommand{\ethm}{\end{Theorem}}
\newcommand{\ble}{\begin{Lemma}}
\newcommand{\ele}{\end{Lemma}}
\newcommand{\bprop}{\begin{Proposition}}
\newcommand{\eprop}{\end{Proposition}}
\newcommand{\bcor}{\begin{Corollary}}
\newcommand{\ecor}{\end{Corollary}}
\begin{document}
\title{Compactness of the automorphism group of a topological parallelism on real projective 3-space}

\author{Dieter Betten and Rainer L\"owen}

\maketitle
\centerline{\it In memoriam Karl Strambach}
\thispagestyle{empty}

\begin{abstract}
We conjecture that the automorphism group of a topological parallelism on real projective 3-space is compact. 
We prove that at least the identity component of this group is, indeed, compact. 
 
MSC 2000: 51H10, 51A15, 51M30 
\end{abstract}

\section{Introduction}

The notion of a topological parallelism on real projective 3-space $\mathop {\rm PG}(3,\Bbb R)$ generalizes the classical example, 
the Clifford parallelism. 
Such a parallelism $\Pi$ may be considered as a set of spreads such that every line belongs to exactly one of them and some continuity property holds.
Many examples of non-classical topological parallelisms have been constructed in a series of papers by Betten and Riesinger, see, e.g., \cite{gl} 
and references given therein. The group $\Phi = \mathop{\rm Aut} \Pi$ of automorphisms of a topological parallelism is a closed subgroup of 
$\mathop {\rm PGL}(4,\Bbb R)$. As is well known, the automorphism group of the Clifford parallelism is the 6-dimensional group 
$\mathop{\rm PSO}(4,\Bbb R) \cong \mathop{\rm SO}(3,\Bbb R) \times \mathop{\rm SO}(3,\Bbb R)$. Betten and Riesinger \cite{coll} proved that no 
other topological parallelism has a group of dimension $\dim \Phi \ge 5$. 
Using the result of the present paper, this has been improved by L\"owen \cite{dim4}: 
in fact, the Clifford parallelism is characterized by $\dim \Phi \ge 4$.
Examples of parallelisms with 1-, 2- or 3-dimensional automorphism groups are known, see \cite{reg}, \cite{gl}, \cite{nonreg}.
Based on the intuition that a parallelism imposes some rigidity akin to an elliptic polarity, we state:\\

\bf Conjecture. \it
The automorphism group of every topological parallelism on $\mathop {\rm PG}(3,\Bbb R)$ is compact.
\rm
\\

This is supported by the observation that the largest possible group, that of the Clifford parallelism, is compact. 
Betten \cite{out} determined the possibilities for groups of dimension $\ge 3$ and found only compact ones. 
Here we shall prove that the assertion of the conjecture is true, at least, for the identity component. This implies the result of \cite{out}.

\bthm\label{main} 
Let $\Pi$ be a topological parallelism on $\mathop {\rm PG}(3,\Bbb R)$ with (full) automorphism group $\Phi$. Then the identity component 
$\Phi^1$ is compact and hence is conjugate to a (closed, connected) subgroup of $\mathop{\rm PSO}(4,\Bbb R)$.
\ethm

The closed, connected subgroups of $\mathop{\rm PSO}(4,\Bbb R)$ are easily enumerated using the fact that $\mathop{\rm SO}(3,\Bbb R)$ 
has only one class of nontrivial proper closed, connected subgroups, represented by the torus group $T = \mathop{\rm SO}(2,\Bbb R)$. 
Up to isomorphism, the nontrivial proper closed, connected subgroups of $\mathop{\rm PSO}(4,\Bbb R)$ are the tori $T$ and $T \times T$, 
$\mathop{\rm SO}(3,\Bbb R)$, and $\mathop{\rm SO}(3,\Bbb R)\times T$. These isomorphism classes split into 3, 1, 3, and 2 conjugacy 
classes, respectively.

We would have liked to produce a proof of Theorem\ref{main} by direct arguments. However, the dynamics of non-compact closed groups 
on real projective 3-space appears to be too diverse to allow for a uniform argument that applies in all situations. Thus, we were forced to 
prove the theorem by examining the conjugacy classes of one-parameter groups one by one. This suffices for the proof because every closed, non-compact
linear group contains a closed, non-compact one-parameter group.

Regarding the status of the conjecture, let us mention that a non-compact group $\Phi$ contains some unbounded cyclic subgroup. 
In most of the possible cases, such a subgroup can be treated (and excluded) by the same methods that we shall use for one-parameter groups. 
There are, however, exceptions; for instance, a cyclic subgroup of a one-parameter group of type (b2) as defined below. So the conjecture 
in full strength remains unproved.
\\

We consider real projective 3-space $\mathop {\rm PG}(3,\Bbb R)$ with its point space $P$ and line space $\CL$. 
Elements of $P$ and $\CL$ are the 1- and 2-dimensional subspaces of the vector space $\Bbb R^4$, respectively. 
The pencil of all lines passing through a point $p$ is denoted $\CL_p$.
The topologies of $P$ and $\CL$ are  the usual topologies of the Grassmann manifolds. See, e.g.,  \cite{Kuehne} for the
continuity properties of geometric operations in this topological projective space. We shall mainly use the facts that $P$ and $\CL$ 
are compact and that the line $L = p \vee q$ joining 
two distinct points $p,q$ depends continuously on the pair $(p,q)$.

\begin{Definition} \rm A \it Parallelism \rm  $\Pi$ on $\mathop{\rm PG}(3,\Bbb R)$ is an equivalence relation on $\CL$ such that for 
every pair $(p,L)\in P \times \CL$ there is a unique line, denoted $\Pi(p,L)$, which passes through $p$ and is parallel (that is, equivalent) to $L$.
The parallelism is said to be \it topological \rm if the map $\Pi : P \times \CL \to \CL$ defined in this way is continuous.
\end{Definition}
 
Using the compactness of $\CL$, it es easily seen that
the continuity condition for $\Pi$ is equivalent to compactness of the parallel relation $\Pi \subseteq \CL \times \CL$. 
In particular, the parallel classes of a topological parallelism are compact spreads and are homeomorphic to the 2-sphere $\Bbb S_2$. 
Every such spread gives rise to a topological affine translation plane with point set $\Bbb R^4$, in which the spread is the set of lines 
through the origin. Here the elements of the spread are considered as 2-dimensional subspaces of $\Bbb R^4$. See \cite{CPP}, Section 64 for details. 
We shall never consider parallelisms that are not topological. Therefore, we shall usually omit the word `topological'.\\

The \it automorphism group \rm $\Phi = \mathop{\rm Aut }\Pi$ is defined as the group of all collineations of $\mathop{\rm PG}(3,\Bbb R)$ preserving 
the equivalence relation $\Pi$. 
We remark that it suffices to interpret `preserving $\Pi$' as the condition that images of parallel lines are again parallel. Then it 
cannot happen that non-parallel lines become parallel after application of the map. Indeed, if a collineation $\gamma$ maps a parallel class
${\cal C} \in \Pi$ into some other parallel class ${\cal C}'$, then $\cal C = C'$, because both sets are spreads. This little observation, 
which is missing in \cite{coll}, ensures that the inverse of an automorphism is again an automorphism and allows us to speak of the automorphism group.
It follows from the compactness of $\Pi$ that $\mathop{\rm Aut }\Pi$ is a closed Lie subgroup of $\mathop {\rm PGL}(4,\Bbb R)$, 
compare \cite{coll}, 4.3.


\section{One-parameter groups: the strategy of proof}\label{EPG}

A non-compact Lie group contains a one-parameter group which is closed and non-compact. This follows, e.g., from the theorem of Mal'cev-Iwasawa, 
\cite{iwa}, Theorem 6. In order to prove our theorem, we shall therefore examine all one-parameter 
subgroups of $\mathop {\rm PGL}(4,\Bbb R)$ and see whether
they can leave a parallelism invariant. It suffices to consider one group from every conjugacy class. Quadratic 
matrices are classified up to conjugacy by their Jordan normal form. The matrix exponential commutes with the conjugation operation. Hence, in order to 
obtain a list of cases to be considered we may use the Jordan normal 
form of a generating $4\times 4$-matrix $A$ for the one-parameter group 
   $$\Gamma = \Gamma_A = \{\exp tA \vert \ t \in \Bbb R\}.$$
Moreover, since only the projective action of $G$ is relevant to us, 
we may replace a generator matrix $A$ by $A + b E$, where $b \in \Bbb R$ and $E$ denotes the unit matrix. This leads to the following list of
cases, where we write $\gamma_t = \exp tA$.\\

\bf Case (a1):  $A$  \it has two complex eigenvalues $\lambda, \mu \in \Bbb C \setminus \Bbb R$ with two Jordan blocks. \rm

Up to addition of a scalar matrix, we have that $\lambda = ai, \mu = b+ci$ with $0\notin \{a,c\}$. Writing $R_\alpha$ for the 
matrix of the rotation of $\Bbb R^2$ through  an angle $\alpha$, we have

   $$\gamma_t = \left(
      \begin{matrix} R_{at} & \cr
                      & e^{bt}R_{ct} \cr 
                      \end{matrix}
         \right)  .            $$

\bf Case (a2): $A$ \it has  one eigenvalue $\lambda \in \Bbb C \setminus \Bbb R$ with a single Jordan block. \rm             

Again we may assume that $\lambda = ai\ne 0$, and we get

 $$\gamma_t = \left(
      \begin{matrix} R_{at} & tR_{at}\cr
                      & R_{at} \cr 
                      \end{matrix}
         \right)  .            $$
   
\bf Case (b1): \it $A$ has one eigenvalue $\lambda \in \Bbb C \setminus \Bbb R$ and two real ones $b, c \in \Bbb R$, the latter with 
Jordan blocks of size 1. \rm

As before, we assume $\lambda = ai \ne 0$ and obtain

 $$\gamma_t = \left(
      \begin{matrix} R_{at} & &\cr
                      & e^{bt} & \cr
                     && e^{ct} \cr 
                      \end{matrix}
         \right)  .            $$
         
\bf Case (b2): \it $A$ has one eigenvalue $\lambda \in \Bbb C \setminus \Bbb R$ and one real eigenvalue $b$, the 
latter with a Jordan block of size 2.    \rm

Taking $\lambda = ai \ne 0$, we get

     $$\gamma_t = \left(
      \begin{matrix} R_{at} & &\cr
                      & e^{bt} & te^{bt} \cr
                     && e^{bt} \cr 
                      \end{matrix}
         \right)  .            $$ 
         
\bf Case (c1): \it $A$ has one real eigenvalue $a$ with a single Jordan block of size 4. \rm         

We may assume that $a = 0$, and we find

     $$\gamma_t = \left(
      \begin{matrix} 1 & t & {1 \over 2} t^2 & {1\over 6} t^3\cr
                      & 1 & t & {1 \over 2}t^2 \cr
                     && 1 & t \cr 
                     &&& 1 \cr
                      \end{matrix}
         \right)  .            $$
         
\bf Case (c2): \it $A$ has two real eigenvalues $a,b$ with Jordan blocks of size 3 and 1, respectively. \rm

We may assume that $a = 0$, and we get 

 $$\gamma_t = \left(
      \begin{matrix} 1 & t & {1 \over 2} t^2 & \cr
                      & 1 & t &  \cr
                     && 1 &  \cr 
                     &&& e^{bt} \cr
                      \end{matrix}
         \right)  .            $$

\bf Case (c3): \it $A$ has two real eigenvalues $a,b$, both with one Jordan block of size 2. \rm    

We may assume that $b = 0$, and then

      $$\gamma_t = \left(
      \begin{matrix} e^{at} & te^{at} & & \cr
                      & e^{at} &  &  \cr
                     && 1 & t \cr 
                     &&  & 1 \cr
                      \end{matrix}
         \right)  .            $$
       
\bf Case (c4): \it $A$ has three real eigenvalues $a,b,c$ with Jordan blocks of size 2, 1, 1, respectively. \rm

Taking $a =0$, we get

 $$\gamma_t = \left(
      \begin{matrix} 1 & t & & \cr
                      &1 &  &  \cr
                     && e^{bt} &  \cr 
                     &&  & e^{ct} \cr
                      \end{matrix}
         \right)  .            $$

\bf Case (c5): \it $A$ is diagonal with eigenvalues $a,b,c,d$. \rm

For $a = 0$, we have

       $$\gamma_t = \left(
      \begin{matrix} 1 &  & & \cr
                      &e^{bt} &  &  \cr
                     && e^{ct} &  \cr 
                     &&  & e^{dt} \cr
                      \end{matrix}
         \right)  .            $$


\section{Eliminating non-compact candidates}\label{elim}

In this section, we shall use the list from the preceding section in order to 
show that $\mathop{\rm Aut} \Phi$ cannot contain any one-parameter group $\Gamma$ that is both 
closed and non-compact. By the theorem of Malcev-Iwasawa \cite{iwa}, this will prove Theorem \ref{main}.

We denote the standard basis vectors of $\Bbb R^4$ by $e_i$, $1 \le i \le 4$. We remind the reader 
that a sequence of vector subspaces $X_n\le \Bbb R^4$ 
of constant dimension $d$ converges to a subspace $X$ in the Grassmann manifold if and only if $X_n$ and $X$ possess bases $x_n^i$, $1 \le i \le d$ and 
$x^i$, $1 \le i \le d$, respectively, such that $x_n^i \to x^i$ as $n \to \infty$.  

\bprop\label{a1}
The group $\Phi$ does not contain a subgroup $\Gamma$ of type (a1) which is closed and non-compact.
\eprop

\bpf Suppose first that $b = 0$. If $c$ is a rational multiple of $a$, then the one-parameter group $\Gamma$ is compact. 
If $c/a$ is irrational, then $\Gamma$ is not closed. In fact, the closure $\overline \Gamma = T$ is a 2-torus and 
is contained in the closed group $\Phi$. 
So the case $b=0$ causes no worries.

Next suppose that $b\ne 0$; reparametrizing, if necessary, we may assume that $b > 0$. 
If $a = c$, let $t_0 = 2\pi/a$. Then $\gamma_{t_0}$ generates an infinite cyclic subgroup and has diagonal shape. 
This possibility can be excluded in the same way as a diagonal one-parameter group of type (c5); see \ref{c5} below. So we have $a \ne c$.

We consider the dynamics of $\Gamma$ acting on $\CL$. Let $H$ be a line not meeting $K = \langle e_1, e_2 \rangle$. 
The images $H^{\gamma_t}$ converge to $L = \langle e_3, e_4\rangle$ as $t \to \infty$. On the other hand, if $M \ne K$ is a line 
meeting $K$ in a point and if the images $M^{\gamma_t}$ accumulate at a line $N$, then $N$ intersects both $K$ and $L$. 
Now consider a parallel class $\CC$ not containing $K$. 
Then $\CC$ contains lines of both types, a line $H$ not meeting $K$ and a 
line $M$ meeting $K$ in a point. (Note that the lines in $\CC$ meeting $K$ form a circle in the 2-sphere $\CC$.)
The $\gamma_t$-images of $H$ and $M$ are always parallel. The images of $H$ converge to $L$ and those of 
$M$ converge to a line $N\ne L$ meeting $L$ in a point. As $\Pi$ is a topological parallelism, $L$ is parallel to $N$, a contradiction.
\phantom{text}
\epf

\bprop\label{a2}
The group $\Phi$ does not contain a subgroup  $\Gamma$ of type (a2).
\eprop

\bpf Setting $t_0 = 2\pi/a$, we obtain a map $\gamma_{t_0}$ generating a cyclic subgroup of a group of type (c3) with $a = 0$. This
cyclic group can be excluded in the same manner as the one-parameter group of type (c3); see \ref{c3special} below. 
\epf

\bprop\label{b1}
The group $\Phi$ does not contain a subgroup  $\Gamma$ of type (b1) which is closed and non-compact.
\eprop

\bpf
If $b = c = 0$, then $\Gamma$ is compact. In the remaining case, the element $\gamma_{t_0}$ with $t_0 = 2\pi/a$ generates an infinite cyclic group of 
diagonal shape which can be excluded by the methods used for type (c5), see \ref{c5} below.
\epf

\bprop\label{b2}
The group $\Phi$ does not contain a subgroup  $\Gamma$ of type (b2).
\eprop

\bpf 
Again we specialize to multiples of $t_0 = 2\pi/a$ and appeal to \ref{c4} below.
\epf

\bprop\label{c1}
The group $\Phi$ does not contain a subgroup  $\Gamma$ of type (c1).
\eprop

Note first that a group $\Gamma$ of type (c1) would fix the line $L = \langle e_1,e_2\rangle$ and hence would also fix the parallel class $\CC$ 
of $L$. Thus $\Gamma$ acts on a translation plane. Such actions are not unusual, however; see, for example, \cite{CPP} 73.10. The proof of 
Proposition \ref{c1} exploits subtle properties of the dynamics and was the hardest one to find. It requires a few preparations. 

We denote the point $\langle e_1 \rangle$ by $p$ and let $\CL_H$ denote the line set of the hyperplane $H = \{x \in \Bbb R^4 \vert \ x_4 = 0\}$.
Finally, $P_H \subseteq P$ denotes the point set of $H$, i.e., the set of its one-dimensional subspaces.

\ble\label{continuous}
Let $\Gamma$ be of type (c1). 

a) On the set $P \setminus P_H$, the maps $\gamma_n$, $n \in \Bbb N$, converge continuously to the constant map with value $p$. 
This means that for every convergent sequence $u_n \to u$ in this set we have $p = \lim_n u_n^{\gamma_n}$.
Likewise, the maps $\gamma_{-n}$ converge continuously to the constant at $p$.

b) On the set $\CL_p \setminus \CL_H$, the maps $\gamma_n$ converge continuously to 
the constant map with value $L$. 
\ele

\bpf a) If $u_n = \langle x_n \rangle$, we may assume that $x_n \to x \ne 0$, and then $u = \langle x \rangle$. In order to find 
$\lim_n u_n^{\gamma_n}$, one has to normalize $y_n = x_n ^{\gamma_n}$ in some way, e.g., by dividing it by its biggest coordinate. 
The claim is then easily verified. Here it is vital that the fourth coordinate  of $x$ is non-zero.

b) Let $X_n \to X$ in $\CL_p\setminus \CL_H$, where $X_n = \langle e_1, x_n\rangle$. We may assume that $x_n \to x \in P\setminus P_H$ and 
that the fourth coordinate $x_{n4}$ of every $x_n$ equals 1. Then we have

       $$X_n^{\gamma_n} = \langle e_1, {2\over n^2}(y_n - y_{n1}e_1)\rangle,$$

where $y_n = x_n ^{\gamma_n}$, as before.  The vectors ${2\over n^2}(y_n - y_{n1}e_1)$ converge to $e_2$, which finishes the proof. We remark 
that another quick proof of assertion (b) can be given by using Pl\"ucker coordinates.
\epf

Assertion (b) of the Lemma fails to hold in general for lines not passing through $p$. It is the discrepancy between 
Lemma \ref{continuous} (b) and this failure that makes the following proof work. \\

\noindent \it Proof of Proposition \ref{c1}. \rm Assume that the group $\Gamma$ of type (c1) leaves $\Pi$ invariant.
We let $p = \langle e_1 \rangle$ as before and let $q$ and $r$ be the points generated by
$e_3$ and $e_4$, respectively. For $n \in \Bbb N$, we define $s_n$ to be the image of $q$ by the inverse map $\gamma_n^{-1} = \gamma_{-n}$. 
By Lemma \ref{continuous} (a), the sequence $s_n$ of points converges to $p$ as $n \to \infty$.
Hence, the line sequence $M_n = s_n \vee r$ converges to the line $M = p \vee r = \langle e_1, e_4 \rangle$, 
which is equal to its own parallel $\Pi(M,p)$. 
The sequence of image lines $M_n^{\gamma_n} = q \vee r^{\gamma_n}$ converges to $q\vee p = \langle e_1, e_3\rangle$; call this line $K$. 
Then $\Pi(M_n^{\gamma_n},p)$ converges to $\Pi(K,p) = K$. Now $\Pi(M_n^{\gamma_n},p)$ equals $(\Pi(M_n,p))^{\gamma_n}$, and $\Pi(M_n,p)$ converges to
$\Pi(M,p) =M \in \CL_p\setminus \CL_H$. Therefore, 
Lemma \ref{continuous} (b) applies and yields  $\lim_n (\Pi(M_n,p))^{\gamma_n} = L$. This is a contradiction.
\phantom{text} \ok
 
\bprop\label{c2}
The group $\Phi$ does not contain a subgroup  $\Gamma$ of type (c2).
\eprop

\bpf
A group $\Gamma \le \Phi$ of type (c2) would fix the line $L = \langle e_1,e_2\rangle$ and the point $r = \langle e_4 \rangle$. Hence it fixes 
the line $\Pi(L,r)$ as well. This means that $\Bbb R^4$ splits as the sum of two invariant 2-dimensional subspaces. This is not 
compatible with the Jordan block structure. 
\epf

For the remaining cases, we need the following lemma, which is of independent interest. 
Here it is convenient to consider $\Pi$ not as an equivalence relation $\Pi \subseteq \CL \times \CL$, but as the set of equivalence classes. In 
the proof of the lemma we shall specify a topology for this set. 

\ble\label{equival}
Suppose that $\Psi \le \Phi = \mathop{\rm Aut} \Pi$ is a subgroup fixing a point $p$ or a hyperplane $H$ of $\mathop {\rm PG}(3,\Bbb R)$. 
Then the action of $\Psi$ on $\CL_p$ or on $\CL_H$, respectively, is equivalent to the action on $\Pi$, considered as the set of parallel classes.
\ele

\bpf 
The maps $\alpha_p: \CL_p \to \Pi$ and $\beta_H: \CL_H \to \Pi$ sending a line to its parallel class are equivalences. 
For bijectivity of the latter map, we 
use the fact that every parallel class is not only a spread, but also a dual spread, which means that every hyperplane contains 
exacty one of its lines; see \cite{CPP}, 64.10a. We introduce a topology on $\Pi$ by insisting that $\alpha_p$ be a 
homeomorphism for a chosen point $p$. Then $\alpha_q$ is a homeomorphism for every point $q$, because 
$\alpha_q^{-1} \circ \alpha_p$ is the homeomorphism $L \mapsto \Pi(L,q)$ from $\CL_p$ to $\CL_q$.  It is now easy to see that 
$\beta_H^{-1}\circ \alpha_p$ is continuous as well, and hence is a homeomorphism by compactness. \phantom{text}
\epf

Note here that the action of $\Pi$ on the pencil 
$\CL_p$ of lines through a fixed point is the same as the action on the projective plane formed by the one-dimensional subspaces of the 
factor space $\Bbb R^4 / p$. 

The following corollary to Lemma \ref{equival} is the `Lemma on Homotheties' of Betten and Riesinger \cite{reg}, 4.2. 

\bcor\label{homoth}
If $\varphi \in \Phi= \mathop{\rm Aut} \Pi$ fixes a point $p$ or a hyperplane $H$ and acts trivially on  $\CL_p$ or on $H$, respectively, 
then $\varphi$ is the identity.
\ecor

\bpf
If $\varphi$ acts trivially on $\CL_p$, then $\phi$ is a central collineation and also possesses an axis, i.e., some hyperplane  is fixed 
pointwise. So assume now that $\varphi$ acts trivially on $H$. 
Every line $L$ has a parallel $M$ contained in $H$.
If $L\ne M$ then $L$ meets $H$ in a point $p$, and $L = \Pi(M,p)$ is fixed by $\varphi$.
\epf

\bprop\label{c3special}
The group $\Phi$ does not contain a subgroup  $\Gamma$ of type (c3) with eigenvalue $a = 0$.
\eprop

\bpf
We have fixed lines $K = \langle e_1,e_2\rangle$ and  $M = \langle e_1,e_3\rangle$. 
The group $\Gamma$ acts on the translation plane $\CT_\CC$ defined by the spread $\CC \in \Pi$ containing $K$ and fixes $M$. As $M$ intersects $K$ 
nontrivially, it cannot belong to $\CC$ and is not a line of this
translation plane, but rather a Baer subplane of $\CT_\CC$. Now the action of $\Gamma$ on $M$ is trivial, but a group fixing a Baer subplane and acting 
trivially on it can have at most two elements \cite{CPP}, 55.21; note that a Baer subplane contains a quadrangle.
\epf

\bprop\label{c3general}
The group $\Phi$ does not contain any subgroup  $\Gamma$ of type (c3).
\eprop

\bpf
By Proposition \ref{c3special} we may assume that $a \ne 0$. 
Let $p = \langle e_1\rangle$ and $q = \langle e_3\rangle$. The action of $\Gamma$ on the projective plane $\CL_p$ is given by the matrices

 $$\left(
      \begin{matrix} e^{at} & &\cr
                      & 1 & t \cr
                     && 1 \cr 
                      \end{matrix}
         \right)  .            $$

The action on $\CL_q$ is projectively equivalent to 

   $$\left(
      \begin{matrix} 1 & t&\cr
                      & 1 & \cr
                     && e^{-at} \cr 
                      \end{matrix}
         \right )  .            $$                 

By Lemma \ref{equival}, these two actions must be equivalent if $\Pi$ is preserved by $\Gamma$. We look at affine versions on $\Bbb R^2$, given by 
$(x,y) \to e^{-at}(x+ty,y)$ and 
$(x,y) \to e^{at}(x+ty,y)$, respectively. Their behaviour for $t$ tending to $+\infty$ is different. 
On the projective plane, both have two fixed points $u,v$.
For $t \to +\infty$, all non-fixed points are moved towards the same fixed point  in one action, and in the other one all non-fixed points except 
those on a single orbit are moved towards the same fixed point, while the remaining ones tend to the other fixed point. Clearly these 
dynamics are not equivalent. 
\epf

\bprop\label{c4}
The group $\Phi$ does not contain a subgroup $\Gamma$  of type (c4).
\eprop

\bpf 
This is an obvious consequence of Lemma \ref{equival}.
\epf

The last remaining case is removed by the following proposition, which therefore completes the proof of Theorem \ref{main}.

\bprop\label{c5}
The group $\Phi$ does not contain a subgroup  $\Gamma$ of type (c5).
\eprop

\bpf 
If three equal eigenvalues occur, then we may assume that these are equal to 0, and  $\Gamma$  acts trivially on the 
hyperplane defined by the eigenspace.
This contradicts Corollary \ref{homoth}. If two eigenvalues are equal, but not three, suppose that $e_1,e_2,e_3,e_4$ are eigenvectors 
for the eigenvalues $0,0,a,b$, respectively. Then the parallel class $\CC$ of the fixed line $\langle e_1,e_3\rangle$  is invariant, 
and $B = \langle e_1,e_2\rangle$ is a Baer subplane of the translation plane defined by $\CC$. The action of $\Gamma$ on $B$ is trivial; 
as in the proof of Proposition \ref{c3special}, this is a contradiction. 

Now suppose that $A$ has distinct eigenvalues $0=a<b<c<d$, with corresponding eigenvectors $e_1,e_2,e_3,e_4$.
Consider the lines $K = \langle e_2,e_4\rangle$ and $L = \langle e_3, e_4\rangle$ and the hyperplanes $G = \{x \vert\ x_1 =0\}$ 
and $H = \{x \vert\ x_4 =0\}$. There is a parallel $M$ of $K$ that is distinct from $K$ and does not meet the line $G\wedge H$ or the line 
$\langle e_1,e_2\rangle$. 
Then $M$ can be written as 
$M = p\vee q$ where $p \in G\setminus H$ and $q \in H\setminus \langle e_1, e_2\rangle$. Now as $t \to \infty$, the images $p^{\gamma_t}$ 
converge to $\langle e_4\rangle$, while $q^{\gamma_t} \to \langle e_3\rangle$. Hence, $M^{\gamma_t} \to L$. On the other hand, $K$
is invariant under $\Gamma$, hence $M^{\gamma_t} = \Pi(K,p^{\gamma_t})$ converges to 
$\Pi (K, \langle e_4\rangle) = K \ne L$, a contradiction.
\epf

We append the first argument that we found for the case of distinct eigenvalues in the previous proof, because we feel that it sheds some extra light on the situation.\\

Suppose that the diagonal matrix $A$ has distinct eigenvalues $0=a<b<c<d$, with corresponding eigenvectors $e_1,e_2,e_3,e_4$. Then the six lines 
$L_{ij} = \langle e_i,e_j\rangle$, $i<j$, are fixed. There are no other fixed lines. Indeed, if $x\in \Bbb R^4$ has three non-zero coordinates,
then the matrix with columns $x, x^{\gamma_t}, x^{\gamma_{2t}}$ contains a $3 \times 3$-matrix obtained from a Vandermonde matrix 
by multiplying each row by a nonzero scalar. Hence these vectors span a 3-space. 

As a consequence, the lines $L_{ij}$ and $L_{kl}$ are parallel whenever $\{i,j,k,l\} = \{1,2,3,4\}$. Indeed, the parallel 
$\Pi(L_{ij},\langle x_k\rangle)$ 
is fixed and disjoint from $L_{ij}$, so this line must be $L_{kl}$.

The group $\Gamma$ acts on the translation plane defined by the spread (parallel class) $\CC_{23}$ containing the fixed lines $L_{23}$ and $L_{14}$. 
As $\Gamma$ is not complex linear, the plane is not the complex plane and its kernel is $\Bbb R$. We may replace $\Gamma$ 
by a subgroup $\Gamma_1$ of $\mathop {\rm SL}(2,\Bbb R)$ without changing the action on $\mathop {\rm PG}(3,\Bbb R)$. 
Then with respect to our translation plane, 
$\Gamma_1$ is contained in a reduced triangle stabilizer in the sense of \cite{CPP}, p. 455, and H\"ahl's theorem on compression 
groups \cite{CPP}, 81.8 applies. Since $\Gamma_1$ contains no non-trivial compact subgroup, the theorem says that 
$\Gamma_1$ itself is a compression group. 
This means that the orbit of every non-fixed line in $\CC_{23}$ has the two fixed lines $L_{23}$ and $L_{14}$ as its ends. But in fact we see that the
orbit of the point $\langle (1,1,1,1)\rangle$ has  end points $\langle e_1\rangle$ and $\langle e_4\rangle$ 
both belonging to $L_{14}$. This is a contradiction.


\bibliographystyle{plain}

\bigskip
\bigskip
\noindent {Dieter Betten,
Mathematisches Seminar der Universit\"at,
Ludewig-Meyn-Str. 4,
D 24108 Kiel, Germany\\ \\
\noindent Rainer L\"owen, Institut f\"ur Analysis und Algebra,
Technische Universit\"at Braunschweig,
Pockelsstra{\ss}e 14,
D 38106 Braunschweig,
Germany}

\end{document}